\documentclass[11pt]{amsart}
\usepackage{amssymb}
\usepackage{amsmath, amsthm, amsopn, amsfonts}
\usepackage[dvips]{graphicx}
\usepackage{graphpap}
\usepackage[english, francais]{babel}
\newtheorem{theoreme}{Theorem}
\newtheorem{theo}{Theorem}[section]

\newtheorem{proposition}{Proposition}
\newtheorem{lemme}[proposition]{Lemma}

\newtheorem{definition}[proposition]{Definition}
\newtheorem{remarque}[proposition]{Remark}

\numberwithin{equation}{section}
\numberwithin{proposition}{section}

\def\11{{\rm 1~\hspace{-1.4ex}l} }
\def\R{\mathbb R}
\def\C{\mathbb C}
\def\Z{\mathbb Z}
\def\N{\mathbb N}
\def\E{\mathbb E}
\def\T{\mathbb T}

\def\cal{\mathcal}
\newcommand{\cqfd}
{%
\mbox{}%
\nolinebreak%
\hfill%
\rule{2mm}{2mm}%
\medbreak%
\par%
}
\begin{document}
\selectlanguage{english}
\title[Estimates for some random series]
{$L^p$ properties for Gaussian random series }
\author{Antoine Ayache}
\address{Laboratoire Paul Painlev\'e, B\^at. M2, Universit\'e Lille 1, 59 655
Villeneuve d'Ascq Cedex, France}
\email{antoine.ayache@math.univ-lille1.fr}
\author{Nikolay Tzvetkov}
\address{
Laboratoire Paul Painlev\'e, B\^at. M2, Universit\'e Lille 1, 59 655
Villeneuve d'Ascq Cedex, France}
\email{nikolay.tzvetkov@math.univ-lille1.fr}
\begin{abstract}
Let $c=(c_n)_{n\in\N^\star}$ be an arbitrary sequence of $l^2(\N^{\star})$ and
let $F_c (\omega)$ be a random series of the type
$$
F_c (\omega)=\sum_{n\in\N^\star}g_n (\omega) c_n e_n ,
$$
where $(g_n)_{n\in\N^*}$ is a sequence of independent ${\cal N}_{\C}(0,1)$
Gaussian random variables and $(e_n)_{n\in\N^\star}$ an orthonormal basis of 
$L^2(Y,{\cal M},\mu)$ (the finite measure space $(Y,{\cal M},\mu)$ being arbitrary).
By using the equivalence of Gaussian moments and an integrability theorem due
to Fernique, we show that a necessary and sufficient condition for $F_c
(\omega)$ to belong to $L^p(Y,{\cal M},\mu)$, $p\in [2,\infty)$ for any $c\in
l^2 (\N^\star)$ almost surely, is that
$\sup_{n\in\N^\star}\|e_n\|_{L^p(Y,{\cal M},\mu)}<\infty$. One of the main
motivations behind this result is the construction of a nontrivial Gibbs
measure invariant under the flow of the cubic defocusing nonlinear
Schr\"odinger equation posed on the open unit disc of $\R^2$.
\end{abstract}

\maketitle
\selectlanguage{english}
%
\section{Introduction and statement of the results}
\label{sec:introduction}

The $d$-dimensional torus ($d\geq 1$) is denoted by
$\T^d\equiv\R^d/(2\pi\Z)^d$. Recall that the exponential functions $\exp(in x)$, $n\in\Z^d$ form a basis of
$L^2(\T^d)$. Let us consider the random multivariate Fourier series
\begin{equation}\label{1}
F_c(\omega,x)=\sum_{n\in\Z^d}g_{n}(\omega)c_{n}e^{inx},
\end{equation}
where $(c_n)_{n\in\Z^d}$ is a non-random sequence of $l^2(\Z^d)$ and $(g_n)_{n\in\Z^d}$ is a sequence of independent ${\mathcal N}_{\C}(0,1)$
Gaussian random variables defined on 
$(\Omega,{\mathcal A},P)$, the underlying probability space. We refer to the
Appendix for some clarifications concerning the series of the type (\ref{1}). 
Any complex-valued centered random variable $X$ whose
probability distribution function (p.d.f.) equals to 
$
\frac{1}{2\pi\sigma^{2}}e^{-|z|^{2}/2\sigma^2}
$
is called an ${\mathcal N}_{\C}(0,\sigma^2)$ Gaussian random variable \footnote{When
$\sigma=0$ then $X\equiv 0$ and its p.d.f. does not exist. However, for convenience, we suppose
that $0$ is a Gaussian random variable}. It is worth noticing that, contrarily
to real-valued Gaussian random variables, the variance of $X$, i.e. the
quantity $\E(|X|^2)$, is not equal to $\sigma^2$ but $2\sigma^2$. Also observe that $X$ can be written as 
$$
X=X_1+iX_2
$$
where $X_1$ and $X_2$ are independent real-valued, centered Gaussian variables
with variance $\sigma^2$. By using a classical result of Paley and Zygmund
one can prove that 
\begin{equation}\label{Lp-finite}
\|F_c(\omega,\cdot)\|_{L^p(\T^d)}<\infty,
\end{equation}
almost surely (a.s.), for any $p\in [2,\infty)$.
The argument to prove (\ref{Lp-finite}) based on the Khinchin inequality uses
that $|\exp(inx)|=1$. 
As we will show in this article, it is however only needed to know that the $L^p(\T^d)$ norms of the 
functions $\exp(in x)$, $n\in\Z^d$ are bounded, uniformly in $n$. 
This  remarkable $L^p$ boundedness of the $L^2$ normalized is a very
particular property of the trigonometric system.
Therefore one can wonder whether the property (\ref{Lp-finite}) remains true
when $\Z^d$ is replaced by $\N^{\star}$, the torus $\T^d$ equipped with the canonical
Lebesgue measure by a finite measure space $(Y,{\cal M},\mu)$ and
the exponential functions by an orthonormal basis $(e_n)_{n\in\N^{\star}}$
of the Hilbert space $L^2(Y,{\cal M},\mu)$.
It turns out that the crucial point  for (\ref{Lp-finite}) to hold 
{\it for every sequence} $c\in l^2(\N^{\star})$ is again the uniform $L^p$ boundedness
of $e_n$. Here is the precise statement.
\begin{theoreme}\label{thm1}
Let us fix $p\in [2,\infty)$ and a  measure space $(Y,{\cal M},\mu)$ with $\mu(Y)<\infty$.
Let $(e_n)_{n\in\N^{\star}}$ be an orthonormal basis of $L^2(Y,{\cal M},\mu)$.
For $c=(c_n)\in l^2(\N^{\star})$, we consider the random series
\begin{equation}\label{2bis}
F_c(\omega)=\sum_{n\in\N^{\star}}g_{n}(\omega)c_{n}e_{n}\,.
\end{equation}
where $(g_n)_{n\in\N^{\star}}$ is a sequence of independent ${\mathcal N}_{\C}(0,1)$
Gaussian random variables.
Then the following two statements are equivalent :
\begin{itemize}
\item[(a)] For every sequence $c\in L^2 (\N^{\star})$,
$$
\|F_c(\omega)\|_{L^p(Y,{\cal M},\mu)}<\infty\mbox{ a.s.}
$$
\item[(b)] One has
$$
\sup_{n\in\N^{\star}}\|e_{n}\|_{L^p(Y,{\cal M},\mu)}<\infty\,.
$$
\end{itemize}
\end{theoreme}
In fact, Theorem~\ref{thm1} is a straightforward consequence of the following
two results.
\begin{theoreme}
\label{Thm:finite}
Under the assumptions of Theorem~\ref{thm1}, 
\begin{equation}\label{finite-serie-ce}
\sum_{n\in\N^{\star}}|c_n|^2 \|e_{n}\|_{L^p (Y,{\cal M},\mu)}^2<\infty
\end{equation}
implies
$
\|F_c (\omega)\|_{L^p (Y,{\cal M},\mu)}<\infty \mbox{ a.s.}
$
\end{theoreme}
\begin{theoreme}\label{Thm:bezkrai}
Under the assumptions of Theorem~\ref{thm1}, 
\begin{equation}\label{eq:bezkrai}
\sum_{n\in\N^{\star}}|c_n|^p \|e_{n}\|_{L^p (Y,{\cal M},\mu)}^p=\infty
\end{equation}
implies
$
\|F_c (\omega)\|_{L^p (Y,{\cal M},\mu)}=\infty \mbox{ a.s.}
$
\end{theoreme}
\noindent {\bf Remarks}
\begin{itemize}
\item[(a)] Theorems~\ref{Thm:finite} and \ref{Thm:bezkrai} mainly rely on the fact
  that a {\em necessary} and sufficient condition for 
$$
\|F_c(\omega)\|_{L^p(Y,{\cal M},\mu )}<\infty \mbox{ a.s.}
$$
to be satisfied is that
\begin{equation}
\label{CNS-finite-norm-Lp}
\E\Big (\|F_c(\omega)\|_{L^p( Y,{\cal M},\mu  )}^p\Big )=\E\Big
(\int_{Y}|F_c(\omega)|^p d\mu\Big )<\infty.
\end{equation}
It is clear that (\ref{CNS-finite-norm-Lp}) is a sufficient condition, but the
fact that it is also a necessary condition may seem surprising. This is
actually a consequence of {\em the Gaussianity of $F_c$ and of an
  integrability theorem due to Fernique }(see Theorem \ref{thm:Fernique}).
\item[(b)] Another important ingredient of the proofs of Theorems \ref{Thm:finite} and
  \ref{Thm:bezkrai} is the {\em the equivalence of Gaussian moments} (see
  Lemma \ref{lem:gauss-moment}. Indeed, it implies that $\mu$ almost every where,
$$
\E(|F_c (\omega)|^p)=C(p)[\E(|F_c (\omega)|^2) ]^{p/2},
$$
where the constant $C(p)$ only depends on $p$. Then it follows from (\ref{2bis})
that 
$$
\E(|F_c (\omega)|^p)=2^{p/2}C(p)\Big
(\sum_{n\in\N^{\star}}|c_n|^2|e_{n}|^2\Big )^{p/2},
$$
which makes the quantity $\E\Big (\|F_c(\omega)\|_{L^p(Y,{\cal M},\mu )}^p\Big )$
less difficult to handle.
\item[(c)] 
Theorem~\ref{Thm:finite} allows the
construction of a nontrivial Gibbs measure on $L^2(D^2)$, invariant under the
flow of the cubic defocusing nonlinear Schr\"odinger equation posed on the
open unit disc $D^2$ of $\R^2$. This was one of the
main motivations behind this study.
We refer to the last section of the paper for details.
\end{itemize}
Let us now make the following definition.
\begin{definition}
Let $c\in l^2 (\N^\star)$ be fixed. The critical $\big(L^p(Y,{\cal M},\mu), (e_n)_{n\in\N^{\star}}\big)$
exponent of the sequence $c$ is defined as
$$
p_{cr}(c)\equiv \sup\Big(p\geq 2\, :\, \|F_c (\omega)\|_{L^p (Y,{\cal M},\mu)}<\infty \mbox{ a.s.}\Big),
$$
or equivalently (in view of Remark (a)) as
$$
p_{cr}(c)\equiv \inf\Big(p\geq 2\, :\, \|F_c (\omega)\|_{L^p (Y,{\cal
    M},\mu)}=\infty \mbox{ a.s.}\Big).
$$
\end{definition}
The previous theorems allow to obtain a lower and an upper bound of
$p_{cr}(c)$, but are not always sufficient to exactly determine this exponent. 
An interesting problem would be to find a general formula allowing to
compute $p_{cr}(c)$ for any sequence $c\in l^2 (\N^\star)$ and any orthonormal
basis $(e_n)_{n\in\N^\star}$ of $L^2(Y,{\cal M},\mu)$. This problem becomes
less difficult to solve when one restricts to some specific situations. One of
them was in fact the starting point of this work and now we are going to
describe it.

Let $D^d\equiv\big\{x\in\R^d\,:\,|x|<1\big\}$ be the unit ball of $\R^d$. Recall that a radial function on $D^d$ is a complex-valued function whose values only
depend on $r\equiv |x|$ the distance to the origin. The subspace of radial functions of the Lebesgue space $L^p(D^d)$, $p\in [1,\infty)$, can be 
identified with the space $L^p([0,1], {\cal B}([0,1]), r^{d-1}dr)$ where
${\cal B}([0,1])$ is the $\sigma$-algebra of Borel sets of $[0,1]$ and $dr$ is the usual Lebesgue
measure on this interval. The eigenfunctions of the Laplace operator on $L^2([0,1], {\cal B}([0,1]), r^{d-1}dr)$, denoted by $e_{n,d}$,
$n\in\N^{\star}$, are closely related to Bessel functions: For any $n\in\N^{\star}$ and real $r>0$, one has  
\begin{equation}\label{def-I-e} 
e_{n,d}(r)\equiv \beta_{n,d}^{-1}r^{-\frac{d-2}{2}}J_{\frac{d-2}{2}}(z_{n,d}\,r),
\end{equation}
where $J_{\frac{d-2}{2}}$ is the Bessel function of order $\frac{d-2}{2}$,
$(z_{n,d})_{n\in\N^\star}$ is the increasing sequence of its (simple) zeroes
on $(0,\infty)$ and 
\begin{equation}\label{beta}
\beta_{n,d}\equiv \Big(\int_{0}^{1}|J_{\frac{d-2}{2}}(z_{n,d}\,r)|^{2}rdr\Big)^{\frac{1}{2}},
\end{equation}
is a normalization factor. Recall that $J_{\nu}$, the Bessel function of
an arbitrary order $\nu\ge 0$, can be defined as the series 
\begin{equation}\label{def-Bessel}
J_{\nu}(r)=\big(r/2\big)^{\nu}
\sum_{n=0}^{\infty}\frac{(-1)^{n}(r/2)^{2n}}{n!\,\Gamma(n+\nu+1)}\, ,
\end{equation}
where $\Gamma (z)=\int_{0}^{+\infty}t^{z-1}e^{-t}\,dt$ denotes the Gamma
function. It follows from the Sturm-Liouville theory that the eigenfunctions 
$e_{n,d}$, $n\in\N^\star$
form an orthonormal basis of $L^2([0,1], {\cal B}([0,1]),
r^{d-1}dr)$. Thus, the analog of (\ref{1}) 
(the eigenfunctions of the Laplace operator on $\T^d$ being $\exp(inx)$, $n\in\Z^d$)
for radial functions on $D^d$ is given by the
random series
\begin{equation}\label{2}
F_{c}(\omega,r)=\sum_{n\in\N^{\star}}g_{n}(\omega)c_{n}e_{n,d}(r),
\end{equation}
where $(c_n)_{n\in\N^\star}$ is a non-random sequence of $l^2(\N^*)$ and
$(g_n)_{n\in\N^\star}$ is again a sequence of independent ${\mathcal N}_{\C}(0,1)$
Gaussian random variables defined on $(\Omega,{\mathcal A},P)$. Here is our
result concerning the computation of the critical $L^p([0,1],{\cal B}([0,1]),
r^{d-1}dr)$ of series of the type (\ref{2}).

\begin{theoreme}\label{th4}
Let $c=(c_n)_{n\in N^\star}\in l^2(\N^\star)$ be a sequence satisfying the
following property :
There exist two constants $0<\alpha_1\le \alpha_2$ such that for any
$n\in\N^\star$,
\begin{equation}\label{dopuskane}
\frac{\alpha_1}{n}\leq c_n\leq \frac{\alpha_2}{n}\,.
\end{equation}
Then the critical 
$\big(L^p([0,1],{\cal B}([0,1]),r^{d-1}dr),e_{n,d}(r)\big)$ exponent of the corresponding series (\ref{2}) is given by
$$
p_{cr}(c)=\frac{2d}{d-2}\,.
$$
\end{theoreme}
\noindent {\bf Remarks}
\begin{itemize}
\item[(a)] 
If we take $c_n=z_{n,d}^{-1}$ then we obtain the series occurring in the Gibbs
measure construction.
\item[(b)] 
One can also calculate $p_{cr}(c)$ for sequences for which (\ref{dopuskane})
is replaced by
$$
\frac{c_1}{n^{\alpha}}\leq c_n\leq \frac{c_2}{n^{\alpha}}\,,\quad \alpha>\frac{1}{2}\,.
$$
\item[(c)] 
For $d\geq 3$, Theorem~\ref{th4} can not be viewed as a consequence of Theorems \ref{Thm:finite} and
\ref{Thm:bezkrai}. Indeed, the latter theorems only imply that $p_{cr}(c)\in
[\frac{2d}{d-2},\frac{2(d-1)}{d-3}]$. 
Actually, in order to compute the precise value of $p_{cr}(c)$, we need to use a precise description of the point-wise
concentrations of $e_{n,d}$ leading to grow of the $L^p(D^d)$ norms of 
$e_{n,d}$.
We refer to \cite{So1,So2,So3} for results giving bounds on the possible
growth of the $L^p$ norms of $L^2$ normalized
eigenfunctions of the Laplace operator on a compact riemannian (boundary-less)
manifold.
\item[(d)]
A function on $S^d$ (the unit sphere $\R^{d+1}$) is called zonal if its value at an arbitrary point $x$ only depends on the geodesic distance
between $x$ and $x_0$, where $x_0$ is a fixed point (not depending on the
function). Roughly speaking, the zonal spherical harmonics (i.e. the
eigenfunctions of the Laplace-Beltrami operator on $S^d$) and the functions  $e_{n,d}$ share the
same properties (see \cite{Stein,Szego}). The latter theorem can therefore be
extended to this new setting.  
\end{itemize}
{\bf Acknowledgement.} We would like to thank Michel Ledoux and Herv\'e
Queff\'elec for their advice. 

\section{Proofs of the main results}

\subsection{Proof of Theorem~\ref{Thm:finite}}
The following elementary lemma is a classical result. It
will be one of the main ingredients of the proof of Theorems \ref{Thm:finite}
and \ref{Thm:bezkrai}.
Roughly speaking it means that the moments of a centered Gaussian random variable are equivalent
\begin{lemme}
\label{lem:gauss-moment}
For any real $p>0$ there exists a constant $C(p)$ such that for any
$Z$, a complex-valued centered Gaussian random variable,
$$
\E(|Z|^p)=C(p)\E(|Z|^2)^{p/2}.
$$
\end{lemme}
As we have already noticed, the quantity
\begin{equation}
\label{fubini-E-F}
\E(\|F_c\|_{L^p(Y,{\cal M},\mu)}^p)=\E\Big(\int_{Y}|F_c(\omega)|^{p}d\mu\Big)=\int_{Y} \E(|F_c(\omega)|^{p})d\mu,
\end{equation}
where $p\in [2,\infty)$ and $c\in l^2 (\N^\star)$ will play a key role in the
proofs of Theorems \ref{Thm:finite} and \ref{Thm:bezkrai}. So let us first express it in a nice way.
\begin{lemme}
\label{lem:nice}
For any exponent $p\in [2,\infty)$ there exists a constant $d(p)$ such that for any sequence $c=(c_n)_{n\in\N^{\star}}$,
\begin{equation}
\label{nice}
\E(\|F_c\|_{L^p(Y,{\cal M},\mu)}^p)=d(p)\Big \|\sum_{n\in\N^\star} |c_n|^2
|e_{n}|^2\Big \|_{L^{p/2}(Y,{\cal M},\mu)}^{p/2},
\end{equation}
where $d(p)>0$ is a constant only depending on $p$.
\end{lemme}
\noindent {\sc Proof of Lemma \ref{lem:nice}:} It follows from part $(b)$ of
Proposition \ref{prop:Banach-F} that 
$\mu$ almost everywhere $F_c(\omega)$ is an ${\cal N}(0,\sigma^2)$ Gaussian
random variable, satisfying 
$$
\E(|F_c(\omega)|^2)=2\sigma^2 =2\sum_{n\in\N^\star} |c_n|^2 |e_{n}|^2.
$$
Thus, by using Lemma \ref{lem:gauss-moment}, one has for any $p\in [2,\infty)$
and $\mu$ almost everywhere, 
$$
\E (|F_c(\omega)|^p)=2^{p/2}C(p)\Big (\sum_{n\in\N^\star} |c_n|^2
|e_{n}|^2 \Big )^{p/2},
$$
which implies that 
\begin{eqnarray*}
\int_{Y}\E (|F_c(\omega)|^p)d\mu&=&
2^{p/2}C(p)
\int_{Y}\Big (\sum_{n\in\N^\star} |c_n|^2|e_{n}|^2 \Big
)^{p/2}d\mu\\
&=&2^{p/2}C(p)\Big \|\sum_{n\in\N^\star} |c_n|^2|e_{n}|^2\Big\|_{L^{p/2}(Y,{\cal M},\mu)}^{p/2}.
\end{eqnarray*}
This completes the proof of Lemma \ref{lem:nice}.
\cqfd
We are now in position to prove Theorem \ref{Thm:finite}. It follows from
(\ref{nice}), the triangular inequality and
(\ref{finite-serie-ce}) that 
$$
\E(\|F_c\|_{L^p(Y,{\cal M},\mu )}^p)\le d(p) \Big (\sum_{n\in\N^\star}
|c_n|^2\|e_{n}\|_{L^p (Y,{\cal M},\mu)}^2 \Big )^{p/2}<\infty.
$$
Thus we get the theorem.
\cqfd
\subsection{Proof of Theorem~\ref{Thm:bezkrai}}
We will make use of the following lemma.
\begin{lemme}\label{nest}
For every $\alpha\geq 1$, every $N\in\N^\star$, every positive numbers
$a_1$,..., $a_N$ one has
$$
\Big(\sum_{n=1}^{N}a_{n}\Big)^{\alpha}\geq\sum_{n=1}^{N}a_{n}^{\alpha}\,.
$$
\end{lemme}
The following integrability theorem, due to Fernique \cite{F}, will play an
important role in our analysis.
\begin{theo}
\label{thm:Fernique}
\cite{F} Let $(E,{\cal B}(E))$ be a measurable vector space (${\cal B}(E)$
denotes the $\sigma$-algebra of Borel sets of $E$), let $F$ be a
centered Gaussian
random variable with values in $E$ and let $N$ be a pseudo-norm on $E$ (the
only difference between a norm and a pseudo-norm is that a pseudo-norm may
take the value $\infty$). If the
probability $P\{N(F)<\infty\}$ is strictly positive, then there exists a
constant $\epsilon_0>0$, such that for any $0<\epsilon<\epsilon_0$, one has
$
\E[\exp(\epsilon N^2(F))]<\infty.
$
In particular, for all $p<\infty$, $\E(N^p(F))<\infty$.
\end{theo}
Let us now give the proof of Theorem~\ref{Thm:bezkrai}.
Using Lemma~\ref{lem:nice}, Lemma~\ref{nest} and the monotone convergence theorem, we infer that
\begin{eqnarray}\label{contradict1}\nonumber
\E(\|F_c\|_{L^p(Y,{\cal M},\mu)}^p) & = & d(p)\Big \|\sum_{n\in\N^\star} |c_n|^2
|e_{n}|^2\Big \|_{L^{p/2}(Y,{\cal M},\mu)}^{p/2}
\\\nonumber
& = &
\lim_{N\rightarrow\infty}
d(p)\Big \|\sum_{n=1}^{N}|c_n|^2
|e_{n}|^2\Big \|_{L^{p/2}(Y,{\cal M},\mu)}^{p/2}
\\\nonumber
& \geq &
d(p)\lim_{N\rightarrow\infty}
\sum_{n=1}^{N}|c_n|^{p}
\|e_{n}\|_{L^{p}(Y,{\cal M},\mu)}^{p}
\\\nonumber
& = &
d(p)\sum_{n=1}^{\infty}|c_n|^{p}
\|e_{n}\|_{L^{p}(Y,{\cal M},\mu)}^{p}
\\
& = &
\infty.
\end{eqnarray}
Then using Theorem~\ref{thm:Fernique} with $E=L^{2}(Y,{\cal M},\mu)$ and $N$
the $L^{p}(Y,{\cal M},\mu)$ norm (viewed as a pseudo-norm on $L^{2}(Y,{\cal M},\mu)$) allows to complete the proof of
Theorem~\ref{Thm:bezkrai}. Indeed, if we suppose that $\|F_c\|_{L^p(Y,{\cal M},\mu)}$ is not a.s. finite then there exists a set $A$ of positive
probability such that for all $\omega\in A$, $\|F_c(\omega)\|_{L^p(Y,{\cal M},\mu)}<\infty$. Thus by the Fernique theorem
$\E(\|F_c\|_{L^p(Y,{\cal M},\mu)}^p)<\infty$ which contradicts (\ref{contradict1}).
\cqfd
\subsection{Proof of Theorem~\ref{thm1}}

Theorem~\ref{thm1} is in fact a consequence of Theorems \ref{Thm:finite} and
\ref{Thm:bezkrai}. It is clear that  Theorem~\ref{Thm:finite} gives the (b) $=>$ (a)
implication of Theorem~\ref{thm1}. Let us now show that (a) implies
(b). Suppose ad absurdum that $(a)$ is satisfied and
\begin{equation}
\label{nonb-th1}
\sup_{n\in\N^\star}\|e_n\|_{L^p(Y,{\cal M},\mu)}=\infty.
\end{equation}
Then it follows from (\ref{nonb-th1}) that there exists a strictly increasing
subsequence $k\mapsto n_k$, such that for every $k\in\N^\star$,
\begin{equation}
\label{subseq-e}
\|e_{n_k}\|_{L^p(Y,{\cal M},\mu)}\ge 2^k.
\end{equation}
Finally let $\tilde{c}$ be the sequence of $l^2(\N^\star)$ defined for every
$n$, as $\tilde{c}_{n}=2^{-k}$ when $n=n_k$ and $\tilde{c}_{n}=0$ elsewhere. By using (\ref{subseq-e}) one obtains that
$$
\sum_{k\in\N^\star}
|c_{k}|^{p}\|e_{n_k}\|_{L^p(Y,{\cal M},\mu)}^{p}=\infty\,.
$$
Then Theorem \ref{Thm:bezkrai} leads to a contradiction.
\cqfd

\subsection{Useful properties of the Bessel functions $J_{\frac{d-2}{2}}$}
In this section we collect several properties of the  Bessel functions useful
in the proof of Theorem~\ref{th4}. 
Let us first recall three important results concerning the zeroes of
$J_{\frac{d-2}{2}}$ and the behaviour of this function. We refer to \cite{S} for their proofs.
 
\begin{itemize}
\item The following inequalities are satisfied for every $n\in\N^\star$ and
  allow to estimate the zeroes of $J_{\frac{d-2}{2}}$:
\begin{equation}\label{U1}
\alpha_{1,d}\,n\leq z_{n,d}\leq \alpha_{2,d}\,n\,,
\end{equation}
where $0<\alpha_{1,d}\le\alpha_{2,d}$ are two constants.
\item A sharp upper bound of $|J_{\frac{d-2}{2}}(r)|$ near the origin is given by:
\begin{equation}\label{U2}
|J_{\frac{d-2}{2}}(r)|\leq C_{1,d}\,r^{\frac{d-2}{2}},
\end{equation}
where $C_{1,d}$ is a constant.
\item The following equality allows to approximate $J_{\frac{d-2}{2}}(r)$:

\begin{equation}\label{U3}
J_{\frac{d-2}{2}}(r)=\sqrt{\frac{2}{\pi\, r}}\cos\Big(r-(d-1)\pi/4 \Big)+R_{d}(r).
\end{equation}
The bigger is $r$ the better is the approximation. Indeed, the remainder $R_{d}$ satisfies
\begin{equation*}
|R_{d}(r)|\leq C_{2,d}\,r^{-\frac{3}{2}},
\end{equation*}
where $C_{2,d}$ is a constant.
\end{itemize} 
Let us now estimate the normalization factor
$\beta_{n,d}$ introduced in (\ref{beta}).
\begin{lemme}\label{beta_eval}
Let $d\geq 1$.
There exist two constants $0<c_1\le c_2$ such that for every $n\in\N^{\star}$,
\begin{equation}\label{estim-beta}
c_1 n^{-\frac{1}{2}}\leq \beta_{n,d}\leq c_2 n^{-\frac{1}{2}}\,\,.
\end{equation}
\end{lemme}
\noindent {\sc Proof of Lemma \ref{beta_eval}:}
It follows from (\ref{beta}) and the variable change $r\rightarrow z_{n,d}\,r$ that
\begin{equation}\label{express-beta}
\beta_{n,d}^{2}=z_{n,d}^{-2}\int_{0}^{z_{n,d}}|J_{\frac{d-2}{2}}(r)|^{2}rdr\,.
\end{equation}
Next using (\ref{U1}), (\ref{U2}) and (\ref{U3}), we get the inequalities
$$
\beta_{n,d}^{2}\leq
Cn^{-2}\Big(1+\int_{1}^{\alpha_{2,d}n}(r^{-\frac{1}{2}})^{2}rdr\Big)\leq
Cn^{-1}\,
$$
and thus we obtain the upper bound part of (\ref{estim-beta}). Let us now
prove that the lower bound part holds. When $n$ big enough,
(\ref{express-beta}) and (\ref{U1}) entail that
$$
\beta_{n,d}\geq
C n^{-1}\left (\int_{1}^{\alpha_{1,d}n}|J_{\frac{d-2}{2}}(r)|^{2}rdr\right)^{1/2}.
$$
Next using (\ref{U3}) and the triangular inequality we get
$$
\beta_{n,d} \geq Cn^{-1}\left (\int_{1}^{\alpha_{1,d}n}\,\,\frac{\cos^{2}\big(r-(d-1)\pi/4
  \big)}{r}rdr\right)^{1/2}-
Cn^{-1}\left (\int_{1}^{\alpha_{1,d}n}r^{-3}rdr\right)^{1/2}.
$$
Finally, using the equality, $\cos^2 a=\frac{1+\cos 2a}{2}$ for any real $a$, we can prove
that 
$$
\lim_{n\rightarrow\infty}
n^{-1}\int_{1}^{\alpha_{1,d}n}\,\,\cos^{2}\big(r-(d-1)\pi/4\big)dr=\frac{\alpha_{1,d}}{2}.
$$
Therefore, we can conclude that for any $n$ big enough, we have
$\beta_{n,d}\geq Cn^{-1/2}$ where $C$ is a constant. This completes the proof of Lemma~\ref{beta_eval}.
\cqfd
We next evaluate the $L^p$ norms of $e_{n,d}$.
\begin{lemme}\label{lemme2}
Let $p\in [2,\infty]$. For convenience we set, for any $n\in\N^\star$,
$$
\delta(n,p,d)=
\left\{
\begin{array}{ll}
1 & {\rm when }\quad 2\leq p<\frac{2d}{d-1},
\\
(\log (2+n))^{\frac{d-1}{2d}} & {\rm when }\quad p=\frac{2d}{d-1},
\\
n^{-\frac{d}{p}+\frac{d-1}{2}} & {\rm when }\quad p>\frac{2d}{d-1}\,.
\end{array}
\right.
$$
Then there exists a constant $C$ such that the inequality,
$$
\|e_{n,d}\|_{L^p(D^d)}\leq C\,\delta(n,p,d)\,,
$$
holds for every $n\in\N^{\star}$.
\end{lemme}
\noindent {\sc Proof of Lemma \ref{lemme2}:}
Let first study the case where $p=\infty$. It follows from (\ref{def-I-e}) that
$e_{n,d}$ can be written as
\begin{equation}\label{U4}
e_{n,d}(r)=\beta_{n,d}^{-1}\,z_{n,d}^{\frac{d-2}{2}}\,G(z_{n,d}\,r),
\end{equation}
where
\begin{equation}
\label{def-G} 
G(r)\equiv r^{-\frac{d-2}{2}}J_{\frac{d-2}{2}}(r).
\end{equation}
Moreover, Relations
(\ref{U2}) and (\ref{U3}) allow to show that $G$ is bounded. On the other hand
Lemma~\ref{beta_eval} and (\ref{U1}) imply that for every $n\in\N^{\star}$,
\begin{equation}\label{bound-beta-z}
C_1 n^{\frac{d-1}{2}}\le\beta_{n,d}^{-1} z_{n,d}^{\frac{d-2}{2}}\leq C_{2}n^{\frac{d-1}{2}}.
\end{equation}
Thus, we obtain the required bound for $p=\infty$. Let us now study the case where
$p<\infty$. Using (\ref{def-I-e}) and Lemma~\ref{beta_eval} one has that
$$
\|e_{n,d}\|_{L^p(D^d)}^{p}
\leq
C_3 n^{\frac{p}{2}}
\int_{0}^{1}
\Big|
r^{-\frac{d-2}{2}}J_{\frac{d-2}{2}}(z_{n,d}\,r)
\Big|^{p}r^{d-1}dr.
$$
Next the variable change $r\rightarrow z_{n,d}\,r$ and (\ref{U1}) yield
$$
\|e_{n,d}\|_{L^p(D^d)}^{p}\leq
C_4 n^{\frac{p}{2}}
n^{-d+p\frac{d-2}{2}}
\int_{0}^{\alpha_{2,d}n}
\Big|r^{-\frac{d-2}{2}}J_{\frac{d-2}{2}}(r)\Big|^{p}r^{d-1}dr.
$$
Thus by using (\ref{U2}) and (\ref{U3}), we get
$$
\|e_{n,d}\|_{L^p(D^d)}^{p}\leq
C_5
n^{-d+p\frac{d-1}{2}}
\Big(
1
+
\int_{1}^{\alpha_{2,d}n}
\,
|r^{-\frac{d-2}{2}}r^{-\frac{1}{2}}|^{p}r^{d-1}dr
\Big).
$$
Next, the inequality, for every $n\in\N^{\star}$,
$$
n^{-d+p\frac{d-1}{2}}\leq (\delta(n,p,d))^{p},
$$
implies that
$$
\|e_{n,d}\|_{L^p(D^d)}^{p}\leq
C_5(\delta(n,p,d))^{p}+
C_5 n^{-d+p\frac{d-1}{2}}
\int_{1}^{\alpha_{2,d}n}
\,
r^{d-1-p\frac{d-1}{2}}dr
\equiv I+II.
$$
Finally, let us upper bound the quantity $II$. Simple computations allow to
show that: If $p>\frac{2d}{d-1}$ then 
$
II\leq 
Cn^{-d+p\frac{d-1}{2}}=C(\delta(n,p,d))^{p}, 
$
If $p<\frac{2d}{d-1}$ then
$
II\leq 
Cn^{-d+p\frac{d-1}{2}}n^{d-p\frac{d-1}{2}}=C(\delta(n,p,d))^{p}
$ 
and if $p=\frac{2d}{d-1}$ then
$
II\leq C\log(2+n)=C(\delta(n,p,d))^{p}.
$
This completes the proof of Lemma~\ref{lemme2}.
\cqfd

Before concluding this subsection let us give a lower bound of $|e_{n,d}|$
near the origin. This bound can be viewed as a measure of the concentration of
$e_{n,d}$ near the origin.
\begin{lemme}\label{lemme3}
There exist two constants $C>0$ and $\gamma>0$, such that for any integer
$n\ge 1$ and real $r\in [0,1]$ satisfying $rn\leq \gamma$, one has 
$$
|e_{n,d}(r)|\geq Cn^{\frac{d-1}{2}}\,.
$$
\end{lemme}
\noindent {\sc Proof of Lemma \ref{lemme3}:}
In view of (\ref{U4}) and (\ref{bound-beta-z}), it sufficient to show that
there exist $C>0$ and $\gamma>0$ two constants, such that for $rz_{n,d}\leq
\gamma$ one has
\begin{equation}\label{U5}
|G(z_{n,d}r)|\geq C.
\end{equation}
It follows from (\ref{def-Bessel}) and (\ref{def-G}) that $G$ can be written
as
\begin{equation}\label{series-G}
G(r)=2^{-\frac{d-2}{2}}
\sum_{n=0}^{\infty}\frac{(-1)^{n}(r/2)^{2n}}{n!\,\Gamma(n+d/2)},
\end{equation}
which implies that $G$ is continuously differentiable on $\R_+$ and
\begin{equation}\label{G-zero}
G(0)=\frac{1}{2^{\frac{d-2}{2}}\Gamma(d/2)}\neq 0.
\end{equation}
Observe that the continuity of $G'$, the derivative of $G$, entails that
\begin{equation}\label{sup-G'}
\sup_{r\in [0,1]} |G'(r)|<\infty.
\end{equation}
Finally, using (\ref{G-zero}), (\ref{sup-G'}) and the equality
$$
G(z_{n,d}r)=G(0)+z_{n,d}\,r\,\int_{0}^{1}G'(t\,z_{n,d}\,r)dt,
$$
one obtains (\ref{U5}). This completes the proof of Lemma~\ref{lemme3}.
\cqfd

\begin{remarque}
\label{rem:infty-e}
Observe that it follows from Lemma \ref{lemme3} that there exists a constant $C>0$
such that for any $n\in\N^\star$, one has 
$$
\|e_{n,d}\|_{L^p (D)}\ge C n^{-\frac{d}{p}+\frac{d-1}{2}}
$$
which implies, when $p>\frac{2d}{d-1}$, that
$$
\lim_{n\rightarrow\infty}\|e_{n,d}\|_{L^p (D)}=\infty.
$$
\end{remarque} 
\noindent {\sc Proof of Remark \ref{rem:infty-e}:}
By using Lemma \ref{lemme3} one has
\begin{eqnarray*}
\|e_{n,d}\|_{L^p (D)}=\Big(\int_{0}^{1}|e_{n,d}(r)|^p r^{d-1}dr\Big
)^{1/p}&\ge & \Big(\int_{0}^{\gamma/n}|e_{n,d}(r)|^p r^{d-1}dr\Big)^{1/p}\\
&\ge & Cn^{\frac{d-1}{2}}\Big (\int_{0}^{\gamma/n}r^{d-1}dr\Big)^{1/p}=Cn^{-\frac{d}{p}+\frac{d-1}{2}}.
\end{eqnarray*}
This completes the proof of Remark~\ref{rem:infty-e}.
\cqfd
\subsection{Proof of Theorem~\ref{th4}}
Before giving the proof of Theorem~\ref{th4} let us observe that
Theorem~\ref{thm1} and the considerations of the previous section show that a new
phenomenon appears when we consider the random series (\ref{2}) instead of
(\ref{1}) : 
There exist sequences $c=(c_n)_{n\in\N^\star}\in l^2
(\N^\star)$ for which, $\|F_{c}(\omega)\|_{L^p (D^d)}<\infty$
a.s.\,, is no longer satisfied by any exponent $p\in [2,\infty)$, but only up to the
critical exponent $p_0=\frac{2d}{d-1}$. More precisely, the following holds
true for $F_{c}$ defined by (\ref{2}) :
\begin{itemize}
\item[(a)]
Let $p<\frac{2d}{d-1}$. Then for every $c=(c_{n})_{n\in\N^{\star}}\in
l^2(\N^{\star})$, one has 
$$\|F_c(\omega)\|_{L^p (D^d)}<\infty \mbox{ a.s.}$$
\item[(b)]
Let $p>\frac{2d}{d-1}$. Then there exists $c\in l^2(\N^{\star})$
such that 
$$\|F_{c}(\omega)\|_{L^p (D^d)}=\infty \mbox{ a.s.}$$
\end{itemize}
Theorem~\ref{th4} is a consequence of Lemma \ref{lemme2}, Theorem
\ref{Thm:finite}, Proposition~\ref{Thm:infinite} and Remark~\ref{alphaz}.
\begin{proposition}
\label{Thm:infinite}
Let $c=(c_n)_{n\in\N^\star}$ be an arbitrary non-zero sequence of $l^2 (\N^\star)$ and
let $\alpha^* (c)$ be the quantity defined as
\begin{equation}
\label{def-alpha-c}
\alpha^* (c)=\sup\Big\{\alpha\ge 0\,:\,\liminf_{N\rightarrow\infty}\Big (
(N+1)^{-\alpha}\sum_{n=1}^N n^{d-1}|c_n|^2\Big) >0\Big\}.
\end{equation}
Then for any exponent $p\in (\frac{2d}{\alpha^* (c)},\infty)$ one has 
$
\|F_c (\omega)\|_{L^p (D^d)}=\infty \mbox{ a.s.}
$
\end{proposition}
Observe that $\alpha^* (c)$ can be viewed as a measure of the speed of
convergence of the series $\sum_{n=1}^\infty |c_n|^2$: The lower it is the 
quicker is the convergence of the series. Also observe that this quantity 
always belongs to $[0,d-1]$. Indeed, one clearly has
$$
\liminf_{N\rightarrow\infty}\Big (\sum_{n=1}^N n^{d-1}|c_n|^2\Big )\ge
\sum_{n=1}^\infty |c_n|^2 >0.
$$
On the other hand, for any $\alpha >d-1$, one has
$$
\liminf_{N\rightarrow\infty}\Big ((N+1)^{-\alpha}\sum_{n=1}^N
n^{d-1}|c_n|^2\Big )=0,
$$
since 
$$
(N+1)^{-\alpha}\Big (\sum_{n=1}^N n^{d-1}|c_n|^2\Big )\le (N+1)^{-(\alpha-d+1)}
\sum_{n=1}^\infty |c_n|^2.
$$
\begin{remarque}\label{alphaz}
Let $c$ be a sequence satisfying (\ref{dopuskane}) then one can directly check
that $\alpha^*(c)=d-2$.
\end{remarque}
\noindent {\sc Proof of Proposition~\ref{Thm:infinite}:} The theorem is clearly
satisfied when $\alpha^* (c)=0$, so we suppose that $\alpha^* (c)>0$. Let $p$
be an arbitrary number of $(\frac{2d}{\alpha^* (c)},\infty)$ and let $\delta$
be a number of $(0,\alpha^*(c))$ verifying
\begin{equation}
\label{ineg-delta-p}
p\ge \frac{2d}{\delta}.
\end{equation}
It follows from (\ref{def-alpha-c}) that there exist $N_0\ge 2$ and a constant
$C_1>0$ such that the inequality 
\begin{equation}
\label{minor-cn}
\sum_{n=1}^N n^{d-1}|c_n|^2\ge C_1 (N+1)^{\delta},
\end{equation}
holds for every $N\ge N_0$. Let us set
$r_0=\min(1,\frac{\gamma}{N_0})$. Observe that any $r\in (0,r_0]$ satisfies
$[\gamma/r]\ge N_0$, where $[\gamma/r]$ denotes the integer part of
$\gamma/r$. 
Thus, putting together Lemma \ref{lemme3} and (\ref{minor-cn}) one obtains, for any
$r\in (0,r_0]$,
\begin{eqnarray}
\label{minor-cn-en}
\nonumber
\sum_{n\in\N^\star} |c_n|^2 |e_{n,d}(r)|^2 &\ge &
\sum_{n=1}^{[\gamma/r]}|c_n|^2 |e_{n,d}(r)|^2\\
\nonumber
&\ge & C\sum_{n=1}^{[\gamma/r]} n^{d-1}|c_n|^2\\
\nonumber
&\ge & C'([\gamma/r]+1)^\delta\\
&\ge & C'' r^{-\delta}.
\end{eqnarray}
Next using (\ref{nice}), (\ref{minor-cn-en}) and (\ref{ineg-delta-p}) one gets
\begin{eqnarray*}
\E(\|F_c\|_{L^p(D^d)}^p) & \ge & d(p) \int_{0}^{1}\Big (\sum_{n\in\N^\star}
|c_n|^2|e_{n,d}(r)|^2\Big )^{p/2}r^{d-1}dr
\\
& \ge &
C_1\int_{0}^{r_0}r^{d-\frac{p\delta}{2}-1}dr
\\
& = &\infty.
\end{eqnarray*}
Next let us fix $\epsilon>0$. Using the last estimates and the inequality
$\exp(\epsilon a^2)\ge C_2 a^p$ for any $a\in\R_{+}$, where $C_2>0$ a constant 
only depending on $\epsilon$ and $p$, one obtains 
$$
\E [\exp (\epsilon\|F_c\|_{L^p (D^d)}]=\infty.
$$
Thus it follows from part $(a)$ of Proposition \ref{prop:Banach-F} and
Theorem \ref{thm:Fernique} that 
$$
\|F_c(\omega,\cdot)\|_{L^p(D^d)}=\infty \mbox{ a.s.}
$$
 This completes the proof of Proposition~\ref{Thm:infinite}.
\cqfd

\section{Application to a invariant measures for the cubic defocusing NLS on the disc }
Let us now describe a consequence of Theorem~\ref{Thm:finite}~ which was one of the
main motivations behind this study.
Theorem~\ref{Thm:finite} allows the
construction of a nontrivial Gibbs measure on $L^2(D^2)$, invariant under the
flow of the cubic defocusing nonlinear Schr\"odinger equation posed on the
open unit disc of $\R^2$, namely the partial differential equation
$$
(i\partial_{t}+\Delta)u-|u|^{2}u=0,
$$
where $\partial_{t}$ is the partial derivative with respect to time, $\Delta$
is the Laplace operator, and  $u(t,x)$ is a complex-valued function defined on
$\R\times D^2$. Using the same method as in \cite{Tz}, one
can show that it sufficient to take the image of the measure 
\begin{equation}\label{invm}
\exp\Big(-\frac{1}{2}\|F_z(\omega)\|_{L^4(D^2)}^{4}\Big)dP(\omega),
\end{equation}
under the function $(\Omega, {\cal A})\rightarrow (L^4 (D^2),{\cal B}(L^4
(D^2)))$,
\begin{equation}\label{map}
\omega\longmapsto F_{z}(\omega,r)\equiv
\sqrt{2}\,\sum_{n\geq 1}\frac{g_n(\omega)}{z_{n,2}}e_{n,2}(r).
\end{equation}
However, in order to show that the latter measure is nontrivial we have to
prove that
\begin{equation}
\label{finite-normL4F}
\|F_{z}(\omega)\|_{L^4 (D^2)}<\infty \mbox{ a.s.}
\end{equation}
Observe that (\ref{finite-normL4F}) cannot be obtained by simply using the
Sobolev embedding,
$$
\|F_z(\omega)\|_{L^4 (D^2)}\le \|F_z(\omega)\|_{H^{1/2}(D^2)}.
$$
Indeed, it follows from (\ref{U1}) and classical properties of i.i.d. centered
Gaussian random variables that
$$
\|F_z(\omega)\|_{H^{1/2}(D^2)}^2\ge C\sum_{n\in\N^\star}
\frac{|g_n(\omega)|^2}{n}=\infty \mbox{ a.s.}
$$
However Theorem~\ref{Thm:finite} yields (\ref{finite-normL4F}).

\section{Appendix}
\label{sec:append}
Let $(Y,{\cal M},\mu)$ be a measure space and $(\Omega,{\cal A},P)$ be a
probability space. For any exponent $p\in [1,\infty]$ we denote by $L^p(Y)$
(resp. $L^p(\Omega)$) the Banach space of complex-valued ${\cal M}$-measurable
functions $f(r)$ defined on $Y$ (resp. of complex-valued random variables $X(\omega)$
defined on $\Omega$) and satisfying $\int_{Y}|f|^p\,d\mu<\infty$ (resp.
$\E(|Y|^p)<\infty$). We denote by $L^p(\Omega\times Y)$ the Banach space of
complex-valued, ${\cal A}\otimes {\cal M}$-measurable functions $Z(\omega,r)$, defined on
$\Omega\times Y$ and satisfying $\E\Big (\int_{Y}|Z|^p\,d\mu\Big)<\infty$. Let 
$\{e_n\}_{n\in\N^{\star}}$ be an orthonormal basis of the Hilbert space
$L^2(Y)$ (we assume that such a basis exists), let $(c_n)_{n\in\N^{\star}}$ be a
non-random sequence of $l^2(\N^{\star})$ i.e. a sequence of complex numbers satisfying 
\begin{equation}\label{cn-l2}
\sum_{n=1}^\infty |c_n|^2<\infty
\end{equation}
and let $(g_n)_{n\in\N^{\star}}$ be a sequence of independent ${\cal N}_{\C}(0,1)$
Gaussian random variables on $(\Omega,{\cal A},P)$. The goal of this section
is to explain why the random series 
\begin{equation}\label{RS-F}
F(\omega,r)=\sum_{n=1}^\infty g_n(\omega)c_n e_n(r),
\end{equation}
is well-defined and to clarify in which sense it converges. Some useful
properties of this random series are also given.
\begin{proposition}
\label{prop:Cauchy-FN}
For any $N\in\N^{\star}$ let $F_N$ be the function $L^2 (\Omega\times Y)$ defined as
\begin{equation}\label{def-FN}
F_N (\omega,r)=\sum_{n=1}^N g_n (\omega)c_n e_n(r).
\end{equation}
Then $(F_N)_{N\in\N^{\star}}$ is a Cauchy sequence in $L^2 (\Omega\times Y)$ and the
series (\ref{RS-F}) is defined as the limit of $(F_N)_{N\in\N^{\star}}$.
\end{proposition}
\noindent {\sc Proof of Proposition \ref{prop:Cauchy-FN}:} By using the fact that the $g_n$'s are independent ${\cal
    N}_{\C}(0,1)$ Gaussian random variables and the fact that
  $\|e_n\|_{L^2(Y)}=1$ for any $n$, one obtains that for every $M<N$,
\begin{eqnarray*}
\|F_M-F_N\|^{2}_{L^2(\Omega\times Y)}
& = &\Big\|\sum_{n=M}^N g_n c_n
e_n\Big\|^{2}_{L^2(\Omega\times Y)}
\\
& = &\sum_{n=M}^N \|g_n\|_{L^2(\Omega)}^2|c_n|^2 
\|e_n\|_{L^2 (Y)}^2
\\
& = &
2\sum_{n=M}^N |c_n|^2.
\end{eqnarray*}
Thus it follows from (\ref{cn-l2}) that
$$
\lim_{M\rightarrow\infty, M<N}\|F_M-F_N\|_{L^2(\Omega\times Y)}=0.
$$
This completes the proof of Proposition \ref{prop:Cauchy-FN}
\cqfd

\begin{proposition}
\label{prop:Banach-F}
The following holds true:
\begin{itemize}
\item[(a)] For $P$-almost all $\omega\in\Omega$, the function
  $F(\omega)$ belongs to $L^2(Y)$ and the $g_n(\omega)c_n$'s are its
  coordinates in the basis $\{e_n\}_{n\in\N^{\star}}$ (on the
  exceptional negligible event, we set $F(\omega,\cdot)\equiv 0$). Moreover, the function
  $F:(\Omega,{\cal A})\rightarrow (L^2(Y),{\cal B}(L^2(Y)))$, $\omega\mapsto
  F(\omega)$, where ${\cal B}(L^2(Y))$ denotes the $\sigma$-algebra of
  Borel sets of $L^2(Y)$, is a centered Gaussian random variable with values
  in the Hilbert space $L^2(Y)$ i.e. $F$ is a measurable function and for any $h\in L^2 (Y)$,
  $\langle F,h\rangle=\int_{Y}F(r)\overline{h(r)}\,d\mu (r)$ is a centered
  complex-valued Gaussian random variable.
\item[(b)] For $\mu$-almost all $r\in Y$ the function $F(\cdot,r)$ is an
  ${\cal N}_{\C}(0,\sigma^2(r))$ random variable, with 
\begin{equation}
\label{sigma-F-r}
 \sigma^2(r)=\sum_{n=1}^\infty |c_n|^2 |e_n (r)|^2.
\end{equation}
\end{itemize}
\end{proposition}

Observe that $\sigma^2(r)<\infty$ for $\mu$-almost all $r$, since one has 
$$2\int_Y \Big (\sum_{n=1}^\infty |c_n|^2 |e_n(r)|^2\Big )\,d\mu (r)=\int_Y \E
\Big (|F(\cdot,r)|^2 \Big )\,d\mu (r)<\infty.
$$

\noindent {\sc Proof of Proposition \ref{prop:Banach-F}:} Let us first prove part $(a)$. As the function
  $(\omega,r)\mapsto F(\omega,r)$ belongs to $L^2 (\Omega\times Y)$, one has
$
\E \Big (\int_Y |F(r)|^2\,d\mu (r)\Big )<\infty
$
and this implies that for $P$-almost all $\omega$, $\int_Y |F(\omega,r)|^2\,d\mu
(r)<\infty$, thus $F(\omega)\in L^2(Y)$. Let us now show that for
$P$-almost all $\omega$ the $g_n(\omega)c_n$'s are the coordinates of
$F(\omega,\cdot)$ in the basis $\{e_n\}_{n\in\N^{\star}}$. In fact it is sufficient to
prove for $P$-almost all $\omega$ one has 
\begin{equation}
\label{gc-l2}
(g_n(\omega)c_n)_{n\in\N^{\star}}\in l^2 (\N^{\star}).
\end{equation}
By using (\ref{cn-l2}) and the fact that the $g_n$'s are independent ${\cal
  N}_{\C}(0,1)$ Gaussian random variables one has that
$$
\E \Big (\sum_{n=1}^\infty |g_n|^2 |c_n|^2\Big )= 2\sum_{n=1}^\infty |c_n|^2
<\infty,
$$ 
which implies that (\ref{gc-l2}) is satisfied a.s. On the exceptional
negligible event ${\cal E}$ where one may have $F(\omega)\notin L^2 (Y)$
or $(g_n (\omega)c_n)_{n\in\N^{\star}}\notin l^2(\N^{\star})$ we set $F(\omega)\equiv 0$.

Let us now prove that the function $F:(\Omega, {\cal
  A})\rightarrow (L^2(Y),{\cal B}(L^2(Y)))$, $\omega\mapsto F(\omega,\cdot)$
is measurable. As $L^2(Y)$ is a separable Hilbert space it is sufficient to
show that for any $u=\sum_{n=1}^\infty a_ne_n\in L^2(Y)$ and any real $\rho$, 
the set
$$
F^{-1}(B(u,\rho))\equiv\Big\{\omega\in\Omega\,:\,\|F(\omega)-u\|_{L^2
  (Y)}<\rho\Big\},
$$
belongs to ${\cal A}$ (note that the sequence $(a_n)_{n\in\N^{\star}}$ belongs to $l^2
(\N^{\star})$ and that $B(u,\rho)$ is the open ball of $L^2(Y)$ of center $u$ a radius
$\rho$). When $0\notin B(u,\rho)$ the set $F^{-1}(B(u,\rho))$ can be expressed as 
\begin{equation}\label{rec-express1}
F^{-1}(B(u,\rho))\equiv\Big\{\omega\in\Omega\setminus {\cal E}\,:\,\sum_{n=1}^\infty |g_n(\omega)c_n-a_n|^2<\rho\Big\},
\end{equation}
and else it can be expressed as 
\begin{equation}\label{rec-express2}
F^{-1}(B(u,\rho))\equiv {\cal E}\cup \Big\{\omega\in\Omega\,:\,\sum_{n=1}^\infty
|g_n(\omega)c_n-a_n|^2<\rho\Big\}.
\end{equation}
By using the fact that 
$$
\sum_{n=1}^\infty
|g_n(\omega)c_n-a_n|^2=\lim_{N\rightarrow\infty}\sum_{n=1}^N
|g_n(\omega)c_n-a_n|^2,\quad {\rm a.s.}
$$ 
one can see that $(\Omega, {\cal A})\rightarrow (\R,{\cal B}(\R))$, $\omega\mapsto
\sum_{n=1}^\infty |g_n(\omega)c_n-a_n|^2$ is a random variable. Thus
(\ref{rec-express1}) and (\ref{rec-express2})
imply that $F^{-1}(B(u,\rho))\in {\cal A}$. Let us now prove that the
measurable function $F:(\Omega, {\cal A})\rightarrow (L^2(Y),{\cal
  B}(L^2(Y)))$, $\omega\mapsto F(\omega)$ is a centered Gaussian random
variable with values in the Hilbert space $L^2(Y)$. Let 
\begin{equation}\label{decomp-h}
h=\sum_{n=1}^\infty
b_n e_n ,
\end{equation}
where $(b_n)_{n\in\N^{\star}}\in l^2 (\N^{\star})$ be an arbitrary function of $L^2 (Y)$. By
using (\ref{decomp-h}) and the fact the $g_n(\omega)c_n$'s are, for $P$-almost
all $\omega$, the coordinates of the function $F(\omega)$ in the basis
$\{e_n\}_{n\in\N^{\star}}$, one has a.s., $\langle F,h\rangle=\sum_{n=1}^\infty g_n c_n
\overline{b_n}$ and this implies that $\langle F,h\rangle$ is a centered complex-valued variable 
(recall that the limit of an almost surely convergent sequence of centered Gaussian
random variables is a centered Gaussian random variable).

Finally, let us prove part $(b)$ of the proposition. It follows from
Proposition \ref{prop:Cauchy-FN} that 
$$\lim_{N\rightarrow\infty}\int_{Y}\E\Big (|F(r)-F_N(r)|^2\Big )\,d\mu(r)=0.$$ 
There exists therefore a sequence $(N_l)_{l\in\N^{\star}}$ such that for $\mu$-almost
all $r$, one has 
$$
\lim_{l\rightarrow\infty}\E\Big (F(r)-F_{N_l}(r)|^2\Big ).
$$
Thus by using the fact that $F_{N_l}(r)$ is an ${\cal
  N}_{\C}(0,\sigma_{N_l}^2(r))$ Gaussian random variable, with
$\sigma_{N_l}^2(r)=\sum_{n=1}^{N_l}|c_n|^2 |e_n(r)|^2$ one obtains part $(b)$
of the proposition (recall that the limit, in the mean square sense, of a
sequence $(\epsilon_n)_{n\in\N^{\star}}$ of ${\cal N}(0,\sigma_{n}^2)$ Gaussian random
variables is an ${\cal N}(0,\sigma^2)$ Gaussian random variable with
$\sigma^2=\lim_{n\rightarrow\infty} \sigma_{n}^2$).
\cqfd

\end{document}